\documentclass[10pt]{article}

\usepackage{amsfonts,amsmath,latexsym,amssymb}    

\usepackage[colorlinks=true,linkcolor=black,citecolor=black,urlcolor=black]{hyperref}

\title{{\smc global geometry under isotropic brownian flows}}

\author{{\smc Sreekar Vadlamani}\thanks{Research supported in part by the Louis and Samuel Seiden Technion Academic Chair}
{\smc \ \
and \ \ Robert J.\ Adler}\thanks{Research supported in part by US-Israel
Binational Science Foundation, Grant 2004064,
 and by Technion VPR Funds.}}

\newtheorem{lemma}{Lemma}[section]
\newtheorem{theorem}{Theorem}[section]

\newcommand{\beq}{\begin{eqnarray}}
\newcommand{\eeq}{\end{eqnarray}}
\newcommand{\beqq}{\begin{eqnarray*}}
\newcommand{\eeqq}{\end{eqnarray*}}
\newcommand{\qed}{\hfill \ensuremath{\Box}}
\newcommand{\be}{\begin{equation}}
\newcommand{\ee}{\end{equation}}
\newcommand{\bea}{\begin{eqnarray}}
\newcommand{\eea}{\end{eqnarray}}
\newcommand{\beaa}{\begin{eqnarray*}}
\newcommand{\eeaa}{\end{eqnarray*}}
\newcommand{\beas}{\begin{eqnarray*}}
\newcommand{\eeas}{\end{eqnarray*}}
\newcommand{\ben}{\begin{enumerate}}
\newcommand{\een}{\end{enumerate}}
\newcommand{\bi}{\begin{itemize}}
\newcommand{\ei}{\end{itemize}}
\newcommand{\ba}{\begin{array}}
\newcommand{\Tube}{{\rm Tube}}
\newcommand{\ea}{\end{array}}
\def\definedas{\stackrel{\Delta}{=}}

\newcommand{\real}{\mathbb{R}}

\newcommand{\lips}{{\cal L}}

\def\LKC{Lipschitz-Killing curvature}

\def\complex{\mathop{\raise .45ex\hbox{${\bf\scriptstyle{|}}$}
     \kern -0.40em {\rm \textstyle{C}}}\nolimits}

\newcommand{\smc}{\scshape}

\def\definedas{\stackrel{\Delta}{=}}

\def\E{{\mathbb{E}}}

\newcommand{\Tr}{{\rm Tr}}

\begin{document}
\baselineskip=12pt

\maketitle


\begin{abstract}

We consider
global geometric properties of a codimension one manifold
embedded in Euclidean space, as it evolves under an isotropic and
volume preserving Brownian flow of diffeomorphisms. In
particular, we obtain expressions describing the expected rate of growth of
the Lipschitz-Killing curvatures, or  intrinsic volumes,
of the manifold under the flow.

These results shed new light on some of the intriguing
growth properties of flows from a global perspective, rather than the
local perspective, on which there is a much larger literature.

\end{abstract}


\section{Introduction}
\label{section.intro}
We are interested in
Brownian flows
 $\Phi_{st},\; 0\le s\le t<\infty$ from $\real^n\to \mathbb{R}^n$,
obtained  by solving the collection of stochastic differential equations
\begin{equation}
\label{equation.flow}
x_t \ =\  \Phi_t(x)\ =\ x + \int_0^t \partial U_s(\Phi_s(x)),
\end{equation}
where we write $\Phi_t$ for $\Phi_{0t}$ when there is no danger of
confusion. Here, $\partial$ denotes the Stratonovich
stochastic differential and $U_t(x)$ is a vector field valued
Brownian motion with smooth spatial covariance structure, on which we shall
have more to say in the subsequent section.
However, we note already that we shall assume
$U$ is such that, with probability one, for each $s\le t$,
\begin{itemize}
\item[(i)] $\Phi_{st}$  is a $C^2$  diffeomorphism.
\item[(ii)]
$\Phi_{st}$ is volume preserving; i.e.\ for any compact
$D\subset\real^n$, $\lambda_n(\Phi_t(D))= \lambda_n(D)$, where
$\lambda_n$ is Lebesgue measure in $\real^n$.
\item[(iii)]
$\Phi_{st}$ is isotropic in the sense of \eqref{isotropy:equation} below.
\end{itemize}
It is standard fare, following from \eqref{equation.flow} and our three
assumptions, that
$\Phi_{ut}\circ\Phi_{su}=\Phi_{st}$, that
$\Phi_{tt}$ is the identity map on $\mathbb{R}^n$,  that
$\Phi_{st}(x)$ and  $\Phi^{-1}_{st}(x)$
are jointly continuous in
$x,s,t$ as are the spatial derivatives
\beq
\label{spatialderiv:equation}
D\Phi_{st}(x)\ \definedas\ \left(\frac{\partial\Phi^i_{st}(x)}{\partial
x^j}\right)_{i.j=1}^n
\eeq
and $D\Phi^{-1}_{st}(x)$,
 and that the `increments'
$\Phi_{s_1t_1},\Phi_{s_2t_2},\ldots,\Phi_{s_nt_n}$ are
independent for all $s_1\le t_1\le s_2\le t_2\le\cdots\le s_n\le t_n$.
For a full study of isotropic Brownian flows, with history and
references, we refer the reader to Kunita's monograph, \cite{Kunbook}.

The study of the evolution of curvature under such flows was pioneered by
LeJan in \cite{LeJan91},
where he established the positive recurrence of the
curvature of a curve moving under an isotropic Brownian flow.
Quite recently, Cranston and LeJan \cite{CrLJ98} followed this with a striking
analysis of the growth of local curvature. Working with isotropic and volume
preserving flows in $\real^n$, they took a codimension one manifold $M$
embedded in $\real^n$, and considered its image under the flow, which we
denote by $ M_t \ \definedas \ \Phi_t(M)$.

Taking a point $x\in M$, they developed an It\^{o} formula
for the symmetric polynomials of the principal
curvatures of $\Phi_t(M)$ at the points $\Phi_t(x)$, including the mean
and  Gaussian curvatures.
In addition, and this will be more important for us, they showed that these
polynomials grow  exponentially in time, with mean rates that are
related to the Lyapunov exponents of the flow.
In simple terms, this means
that the manifold $M_t$, while it may begin at time zero as something
as simple as the unit sphere $S^{n-1}$ in $\real^n$ (which has unit
Gaussian curvature everywhere) it tends to develop sharply
rounded `corners' as time progresses.

A somewhat different set of results can be found in a series of papers
\cite{CrScSt99,CSS} authored by
Cranston, Scheutzow and Steinsaltz.
In particular, the
combined results of
\cite{CrScSt99,CSS,SS2001} show that,
for an isotropic Brownian flow with $n\geq 2$, there are positive constants $c$ and $C$
such that for each
compact and connected set $D\subset\real^n$ with at least two points,
\beq
\label{fastgrowth:equation}
c\ \leq \ \inf_{u\in S^{n-1}} \sup_{x\in D} \liminf_{t\to\infty}
\frac{1}{t} \langle \Phi_t(x),u\rangle \ \leq\  \limsup_{t\to\infty}
\frac{1}{t} \sup_{x\in D} \|\Phi_t(x)\|\leq  C,
\eeq
almost surely if the top Lyapunov exponent is strictly positive,
 and with strictly positive probability, otherwise.

One implication of this result is that while $\Phi_t(M)$
is  homotopically equivalent to $M$, for large $t$ it will `look' quite
different. One way to measure this difference, at a global level, is via
their \LKC s. Since $M$ has dimension $(n-1)$ there are $n$ such
curvatures, $\lips_0(M_t),\lips_1(M_t),\dots,\lips_{n-1}(M_t)$. The
first of these, $\lips_0(M_t)$, is the Euler-Poincar\'e
characteristic of $M_t$, which, because $M$ and $M_t$ are
homotopically equivalent, is the same as that of $M$, and so
independent of $t$. The last of these, $\lips_{n-1}(M_t)$, gives the
$(n-1)$ dimensional surface measure of $M_t$ and most definitely
does change with time, as do all the remaining $\lips_j(M_t)$.
Further information on the geometric r\^oles of the \LKC s is given
in the following section.

In view of the results of Cranston, Scheutzow and Steinsaltz described above,
one would expect that the  $\lips_j(M_t)$, $1\leq j\leq n-1$, would grow
rapidly in time, as parts of the set $M_t$ begin to stretch in various
directions at rate $t$. That this is indeed the case is a consequence of the
following theorem, one of the two main results of this paper.

\begin{theorem}
\label{theorem.mean.lk}
Let $M$ be a smooth codimension one manifold embedded in $\mathbb{R}^n$ and $M_t$
its image under $\Phi_t$, where $\Phi_{st}$ is an
isotropic and volume preserving Brownian  flow of $C^2$ diffeomorphsims
of $\real^n$.
Then, for  $0\leq k\leq n-1$, the expected rate
of growth of the Lipschitz-Killing curvatures is given by
\begin{equation}
\label{equation.meanLK}
\E\left\{{\cal L}_{n-k-1}(M_t)\right\} \ =\
 {\cal L}_{n-k-1}(M)\exp\left(\frac{(n-k-1)(n+1)(k+1)\mu_2t}{2n(n+2)}\right),
\end{equation}
where $\mu_2$ is the second moment of the spectral measure $F$ of
\eqref{spectralmeasure:equation}.
\end{theorem}

The proof of Theorem \ref{theorem.mean.lk} relies on the fact that,
loosely speaking, for $1\le k\le (n-1)$, the $(n-k-1)$-th
Lipschitz-Killing curvature of a manifold can be obtained as an
average, over the manifold, of the $k$-th order symmetric polynomial
of the principal curvatures. However, as we have already noted
above, these have been studied in detail by Cranston and LeJan
\cite{CrLJ98}. Consequently, our proof relies very heavily on their
paper, to the extent that one could consider this paper as an
addendum to theirs. Nevertheless, we believe that the results are of
independent interest, in that they lift the local approach of
\cite{CrLJ98} to a global scenario. Finally the case $k=0$, which
corresponds to the $(n-1)$ dimensional surface measure of the
manifold or equivalently the $(n-1)$-th Lipschitz-Killing curvature,
is a simple consequence of Lemma \ref{theorem.ito.norm}.

The remainder of this paper is organized as follows. In
Section~\ref{section.lkc} we provide the precise definition of
Lipschitz-Killing curvatures and the required geometric background,
followed by the proofs of the main results of the paper in
Section~\ref{section.main}.


\section{Brownian flows}
\label{section.background}

This section is not so much about Brownian flows {\it per se}, for which
we refer you back to the references of the Introduction, but rather
about setting up notation. Since we plan to use the main result of
\cite{CrLJ98} to prove Theorem \ref{theorem.mean.lk}, we shall adopt the
notation of that paper without much explanation. You can find
missing explanations in \cite{CrLJ98}.

The first step is to define the vector field valued Brownian motion
$U$ driving the flow in \eqref{equation.flow}. We take this to be a
zero mean Gaussian process from $\real_+\times\real^n$ to $\real^n$
with covariance structure given by
\begin{equation*}
\E\left\{  U_t^k(x) U_s^l(y)\right\}\ =\ (t\wedge s)C^{kl}(x-y), \;\;\;\; 1\le k,l\le n,
\end{equation*}
where each $C^{kl}$ can be written in the form
\begin{equation}
\label{spectralmeasure:equation}
C^{kl}(z) \ =\ \int_0^\infty \int_{S^{n-1}} e^{i\rho \langle
z,t\rangle} (\delta^k_l -t^k t^l) \, \sigma_{n-1} (dt) F(d\rho),
\end{equation}
for a normalized Lebesgue measure $\sigma_{n-1}$ on $S^{n-1}$
and a non-negative measure $F$ on $\mathbb{R}^{+}$, with $m$-th moment,
$\mu_m$, given by $\mu_m=\int_0^\infty \rho^m\,F(d\rho)$.

The various spatial derivatives of $U$,
implicitly assumed to exist, are denoted by
\begin{eqnarray*}
W^i_j \ = \ \frac{\partial U^i}{\partial x^j},\quad
B^i_{jk} \ =\ \frac{\partial^2 U^i}{\partial x^j \partial x^k}.
\end{eqnarray*}
Writing $\langle\cdot,\cdot\rangle$ for quadratic covariation, it is not hard to  check that
\begin{eqnarray}
\label{equation.quad.cov}
\langle dW^i_j(t,y),dW^k_l(t,y)\rangle & = & \frac{\mu_2}{n(n+2)}
\left[(n+1)\delta^i_k\delta^j_l-\delta^i_j\delta^k_l- \delta^i_l\delta^k_j\right]dt,\\
\langle dB^i_{jk}(t,y), dW^p_q(t,y) \rangle & = & 0,
\end{eqnarray}
for any $1\le i,j,k,l,p,q\le n$, implying $E(\sum W^i_i)^2=0$, hence volume preserving, and
\beqq
\langle\langle dB(u,u),v\rangle,\langle dB(u,u),v\rangle\rangle
= \frac{3\mu_4}{n(n+2)(n+4)}\left[(n+3)\|u\|^4\|v\|^2 -4\langle
u,v\rangle^2\|u\|^2\right]dt,
\eeqq
 for all vectors $u,v \in \mathbb{R}^n$.

This particular choice of the covariance function makes the flow
isotropic in the sense that the spatial covariance matrices
$C(x)=(C^{kl}(x))_{l,k=1}^n$ satisfy
\beq
\label{isotropy:equation}
C(x) \ = \ G^*C(Gx)G
\eeq
for any real orthonormal matrix $G$, as well as making the flow volume preserving.


With the flow defined, we now turn to setting up the notation required
for studying its (differential) geometry. Our basic references for this
are Lee \cite{Lee.Rie,Lee.Smooth} and
Part II of  \cite{GRF}.

We start with a
codimension one Riemannian manifold $M$ embedded in $\mathbb{R}^n$
%
and, for $x\in M$, let $u=\{u_i\}_{i=1}^{(n-1)}$ be an orthonormal
basis of $T_xM$, the tangent space at $x\in M$. (Note that $u$
actually depends on $x$, but we shall not write this explicitly.)
Then by a simple {\it push-forward} argument,
$D\Phi_t(x)u_i=u_i(t)\in T_{x_t}M_t$, where $x_t=\Phi_t(x)$ (cf.\
\eqref{equation.flow}) and $D\Phi_t(x)$ is as defined in
\eqref{spatialderiv:equation}. Furthermore, \beqq du_i(t)\ =
\partial Wu_i(t) = dWu_i(t). \eeqq

Writing $\Pi_t:\ T_{x_t}\mathbb{R}^n \to T_{x_t}M_t$ as the
orthogonal projection onto $T_{x_t}M_t$, and
$\widetilde{\nabla}$ for the canonical connection on the ambient
Euclidean space $\real^n$,
the second fundamental form at
$x_t \in M_t$ is given by
\beqq
S_t(u(t),v(t))\ =\ \left(I-\Pi_t\right)\widetilde{\nabla}_{u(t)}v(t),
\eeqq
for any  $u(t), v(t)\in T_{x_t}M_t$.
Subsequently,
\emph{scalar} second fundamental form is defined as
\beqq
S_{\nu_t}(u(t),v(t))\ =\ \langle S_t(u(t),v(t)),\nu_t\rangle,
\eeqq
where   $\nu_t$ denotes the unit normal vector field and
$\langle\cdot,\cdot\rangle$ now denotes the usual Euclidean inner product
rather than quadratic covariation.

It is standard fare that the scalar
second fundamental form can be used to induce a
linear operator on the exterior algebra $\Lambda^k(T_{x_t}M_t)$ of
alternating covariant tensors on $T_{x_t}M_t$, built over the usual wedge
product,
for each $1\le k\le (n-1)$.
First we define $S^{(k)}_{\nu_t}$ as \beqq
S^{(k)}_{\nu_t}(u_{l_1}(t)\wedge\ldots\wedge
u_{l_k}(t),\;u_{m_1}(t)\wedge\ldots\wedge u_{m_k}(t))\ = \
\sum_{\sigma\in S_k}(-1)^{\eta_\sigma}\prod_{j=1}^k
S_{\nu_t}(u_{l_{\sigma(j)}}(t),u_{m_j}(t)), \eeqq where
$\{u_{l_p}(t)\}, \{u_{m_q}(t)\}\subset \{u_i(t)\}_{i=1}^{(n-1)}$,
$S_k$ is the collection of all permutations of $\{1,\dots,k\}$ and
$\eta_\sigma$ denotes the sign of the permutation $\sigma$. This
gives rise to the linear operator
\beqq
u_{l_1}(t)\wedge\ldots\wedge
u_{l_k}(t)\ \mapsto\ S^{(k)}_{\nu_t}(u_{l_1}(t)\wedge\ldots\wedge
u_{l_k}(t),\cdot), \eeqq
 where
$\{u_i(t)\}_{i=1}^{(n-1)}\subset T_{x_t}M_t$ is a basis of
$T_{x_t}M_t$.

The last and the most important remaining definition is that of
the trace of $S^{(k)}_{\nu_t}$. For this, however, we need
some more notation. For $1\le k\le (n-1)$ define the index set $I_k$ by
\beqq
I_k\ =\ \left\{\vec{m}\in\{1,\ldots,n-1\}^k:\: m_1<m_2<\cdots<m_k\right\}.
\eeqq
Then, for $\vec{l}\in I_k$, define
\beqq |\vec{l}| &=& l_1 +\cdots +l_k,\\
\label{equation.wedge.prod} \alpha_{\vec{l}}(t)&=&u_{l_1}(t)\wedge
\cdots \wedge u_{l_k}(t),
\\
\alpha^{\vec{l}}(t)&=&(-1)^{|\vec{l}| +k} u_1(t)\wedge\cdots\wedge{\widehat u}_{l_1}(t)\wedge\cdots\wedge{\widehat
u}_{l_k}(t)\wedge\cdots\wedge u_{n-1}(t),\\
\alpha(t) &=& u_1(t)\wedge\cdots\wedge u_{n-1}(t),
\eeqq
where the hatted vectors are
understood to be omitted from the wedge product. Now, for
$\vec{l},\vec{m}\in I_k$, define
\begin{equation}
\label{equation.inner.prod}
\langle\alpha_{\vec{l}}(t),\alpha_{\vec{m}}(t)\rangle
\ =\ \det\left(\langle u_{l_i}(t),u_{m_j}(t)\rangle\right),
\end{equation}
and naturally
$\|\alpha(t)\|^2 \ =\ \det\left(\langle u_i(t),u_j(t)\rangle\right)$.

We now have all that we need to define the
all important trace, $\Tr S^{(k)}_{\nu_t}$, as
\begin{equation}
\label{equation.trace}
\Tr S^{(k)}_{\nu_t}\  =\
\frac{\langle
\alpha^{\vec{l}}(t),\alpha^{\vec{m}}(t)\rangle}{ \|\alpha(t)
\|^{2}} \, S^{(k)}_{\nu_t}(\alpha_{\vec{l}}(t), \alpha_{\vec{m}}(t)),
\end{equation}
where the Einstein summation convention is carried over the indices
$\vec{l},\vec{m}\in I_k$. Now for the case $k=0$, which has thus far
remained untouched, we define $Tr S^{(0)}_{\nu_t}=1$.

The temporal development of this trace was studied in detail in
\cite{CrLJ98}, where Cranston and LeJan proved
that $\{\Tr S^{(k)}_{\nu_t}\}_{k=1}^{(n-1)}$ is a
$(n-1)$-dimensional diffusion. This is the precise formulation of the
result we were referring to
in the Introduction when we spoke of the It\^{o} formula of symmetric polynomials of
principal curvatures, which are essentially equivalent to the traces
$\{\Tr S^{(k)}_{\nu_t}\}_{k=1}^{(n-1)}$.


\section{Lipschitz-Killing curvatures}
\label{section.lkc}

There are a number of ways to define \LKC s, but perhaps the easiest
is via their appearance in the so-called tube formulae, which, in
their original form, are due to Weyl \cite{Wey39}. (For more details
and applications see either the monograph of Gray \cite{GRAY} or
Chapter 10 of \cite{GRF}.)

To state the tube formula, let $M$ be a $C^2$,
$(n-1)$-dimensional manifold
embedded in $\mathbb{R}^n$ and endowed with the canonical
Riemannian structure on $\mathbb{R}^n$. The
tube of radius $\rho$ around $M$ is defined as
\beqq
\Tube(M,\rho) \ = \ \left\{x\in\mathbb{R}^n:\; d(x,M)\le\rho \right\},
\eeqq
where
\beqq
d(x,M)\ = \ \inf_{y\in M} \|x-y\|.
\eeqq
Weyl's
tube formula states that there exists a $\rho_c\ge 0$, known as
the critical radius of $M$, such that, for $\rho\le\rho_c$, the
volume of the tube is given by
\beq
\label{LKC:definition}
\lambda_n(\Tube(M,\rho))\ = \ \sum_{j=0}^{n-1}\rho^{n-j}\omega_{n-j}{\cal L}_j(M),
\eeq
where $\omega_j$ is the volume of the $j$-dimensional unit ball and
${\cal L}_j(M)$ is the $j^{th}$-Lipschitz-Killing curvature of $M$.

Writing ${\cal H}_j$ for $j$-dimensional Hausdorff measure, it is
easy to check from \eqref{LKC:definition} that $\lips_{n-1}(M)={\cal
H}_{n-1}(M)$. That is, it is the surface `area' of $M$.
$\lips_{0}(M)$ is the Euler-Poincar\'e characteristic of $M$, and
while the remaining \LKC s have less transparent interpretations, it
is easy to see that they satisfy simple scaling relationships, in
that $\lips_j(\alpha M) = \alpha^j \lips_j(M)$ for all $1\leq j\leq
n-1$, where $\alpha M=\{x\in\real^n\: x=\alpha y\ \text{for some}\
y\in M\}$. Furthermore, despite the fact that defining the $\lips_j$
via \eqref{LKC:definition} involves the embedding of $M$ in
$\real^n$, the ${\cal L}_j(M)$ are actually intrinsic, and so
independent of the embedding space.

While \eqref{LKC:definition} characterizes the $\lips_j(M)$ it does
not generally help one compute them. There are a number of ways in
which to do this, but we choose the following, which is most
appropriate for our purposes. (cf.\ \cite{GRF,GRAY} for further
details and examples)
\beq
{\cal L}_{n-k-1}(M) \ =\
 K_{n,k}\int_{M}\int_{S(\mathbb{R})}\Tr S^{(k)}_{\nu}
 1_{N_xM}(-\nu)\, {\cal H}_0(d\nu){\cal H}_{n-1}(dx),
\label{LK:equation}
\eeq
where
$ K_{n,k}\ =\ \frac{1}{2(\pi)^{(k+1)/2}k!}\Gamma\left(\frac{k+1}{2}\right)$ and
$N_xM$ is the normal cone to $M$ at the point $x$.
Since $M$ has codimension 1 in $\real^n$, each
$N_xM$ contains only the outward normals to $T_xM$ and is of unit dimension.

In general, we shall write $S(\mathbb{R}^n)$ for the unit sphere in
$\real^n$. Thus the
$S(\mathbb{R})$  appearing in \eqref{LK:equation} contains only the
two vectors $+1$ and $-1$ in $\real$ and ${\cal H}_0$ is counting measure,
which makes the integral over $S(\mathbb{R})$ a rather pretentious way of
writing things, as, indeed, was the introduction of the normal cone.
Nevertheless, both will be of use to us later on, when we
discuss possible generalisations of our results.

We, of course, are interested in the temporal evolution of the
$\lips_j(M_t)$ and it follows directly from \eqref{LK:equation} that,
with the notation of the previous section,

\begin{eqnarray}
\label{equation.LK}
{\cal L}_{n-k-1}(M_t) & = & K_{n,k}\int_{M_t}\int_{S(\mathbb{R})} \Tr S^{(k)}_{\nu_t}
1_{N_{x_t}M_t}(-\nu_t)\,{\cal H}_0(d\nu_t){\cal H}_{n-1}(dx_t)\nonumber\\
& = & K_{n,k}\int_{M}\int_{S(\mathbb{R})} \Tr S^{(k)}_{\nu_t} \sqrt{\det(\langle u_i(t),u_j(t)\rangle)}
1_{N_{x_t}M_t}(-\nu_t)\,{\cal H}_0(d\nu_t){\cal H}_{n-1}(dx)\nonumber\\
& = & K_{n,k}\int_M\int_{S(\mathbb{R})} \Tr S^{(k)}_{\nu_t} \|\alpha_t\|
1_{N_{x_t}M_t}(-\nu_t)\,{\cal H}_0(d\nu_t){\cal H}_{n-1}(dx).
\end{eqnarray}



\section{An It\^o formula for ${\cal L}_{n-k-1}(M_t)$}
\label{section.main}

Before commencing a serious stochastic analysis of
\eqref{equation.LK} we recall some of the results and further
notation from LeJan \cite{LeJan85}.

Let $\xi(t)=\xi_1(t)\wedge\ldots\wedge \xi_k(t)$ and $\psi(t)=\psi_1(t)\wedge\ldots\wedge \psi_k(t)$,
where $\{\xi_i(t)\},\{\psi_i(t)\}\subset T_{x_t}M_t$. Then, by
Lemma 3 of \cite{LeJan85},
\begin{eqnarray*}
d\langle\xi(t),\psi(t)\rangle =  \sum_{l,\;j}(\langle\tau_i^j\xi(t),\psi(t)\rangle +\langle\xi(t),\tau_i^j\psi(t)\rangle)\,dW^i_j(t)
\, +\, \frac{k(n-k)\mu_2}{n}\langle\xi(t),\psi(t)\rangle \,dt,
\end{eqnarray*}
where
\beqq
\tau^j_l\xi(t) \ =\  e^j\wedge \left\{\sum_{i=1}^k
  (-1)^{i+1}\langle\xi_i(t),e^l\rangle\xi_1(t)\wedge\cdots\wedge{\hat
    \xi}_i(t)\wedge\cdots\wedge\xi_k(t)\right\},
\eeqq
with $\{e^{k}\}_{k=1}^n$ being the standard basis of $\mathbb{R}^n$, and
$\tau^j_l\psi(t)$ is defined similarly.

It follows immediately from the above that if
$\xi(t)=\psi(t)=\alpha(t) =u_1(t)\wedge\cdots\wedge u_{n-1}(t)$, where $u_i(t)=D\Phi_t(x)u_i$ and $(u_1,\ldots,u_{n-1})$ is an orthonormal
basis of $T_xM$, then
\beqq
d\|\alpha(t)\|^2\ =\ \|\alpha(t)\|^2\Big(2\sum_{i=1}^{n-1}dW^i_i(t)\, +\,
\frac{(n-1)\mu_2}{n}dt\Big).\eeqq

Now we  derive an It\^{o} formula for $\|\alpha(t)\|$
and use it together with Theorem A.2 of \cite{CrLJ98} to
obtain an expression for the It\^{o} derivative of the Lipschitz-Killing
curvatures.

\begin{lemma}
\label{theorem.ito.norm}
Let $M$ be a smooth $(n-1)$-dimensional manifold embedded in $\mathbb{R}^n$ and $M_t$ its image at time $t$ under the stochastic, isotropic, and volume preserving flow $\Phi_t$ described in Section \ref{section.background}. Then,
in the notation of Section \ref{section.background},
\beqq
d\|\alpha(t)\| \ =\  \|\alpha(t)\|\Big(\sum_{i=1}^{n-1}dW^i_i(t)\, +\,\frac{(n-1)(n+1)\mu_2}{2n(n+2)}dt\Big).
\eeqq
\end{lemma}

\textbf{Proof:} Using the standard It\^{o} formula and \eqref{equation.quad.cov} we obtain
\begin{eqnarray*}
d\|\alpha(t)\| & = & d(\|\alpha(t)\|^2)^{\frac{1}{2}}\\
 & = & \frac{1}{2}\|\alpha(t)\|\Big(2\sum_{i=1}^{n-1}dW^i_i(t)
 +\frac{(n-1)\mu_2}{n}dt\Big)
-\ \frac{(n-1)\mu_2}{2n(n+2)}\|\alpha(t)\|dt\\
        & = & \|\alpha(t)\|\Big(\sum_{i=1}^{n-1}dW^i_i(t) +\frac{(n-1)(n+1)\mu_2}{2n(n+2)}dt\Big).
\end{eqnarray*}
\qed

We need just a little more preparation before  we can turn to
 our main result.

Let $\alpha_{\vec{l}}(t)=u_{l_1}\wedge\cdots\wedge u_{l_k}(t)$, be a
$k$-form for $\vec{l}\in I_k$, then define
\begin{equation}
\label{equation.sub.wedge.product}
\alpha_{\vec{l_p}}(t) = (-1)^{p+1}u_{l_1}(t)\wedge\cdots\wedge{\hat u}_{l_p}(t)\wedge\cdots\wedge u_{l_k}(t),
\end{equation}
for $1\le k\le (n-1)$, where $\vec{l}\in I_k$, $\vec{l_p}\in I_{k-1}$ and $1\le p\le k$.

Rewriting the above expression as
\begin{equation}
\label{equation.sub.wedge.product1}
\alpha_{\vec{l_p}}(t) \ =\  (-1)^{p+1}u_{l_1}^{(p)}(t)\wedge\cdots\wedge u_{l_{k-1}}^{(p)}(t),
\end{equation}
defines $u_l^{(p)}$.

Then according to Theorem A.2 of \cite{CrLJ98}, for $1\le k\le n-1$
\begin{eqnarray}
\label{equation.ito.trace}
d\Tr S^{(k)}_{\nu_t} & = & \sum_{i,p}\left[ S^{(k-1)}_{\nu_t}(\alpha_{\vec{l_p}}(t),\alpha_{\vec{m_i}}(t))\langle dB(u_{l_p}(t),u_{m_i}(t)),\nu_t\rangle\right] \frac{\langle\alpha^{\vec{l}}(t),\alpha^{\vec{m}}(t)\rangle}{\|\alpha(t)\|^{2}} \nonumber\\
           &   & \quad +\ \Tr S^{(k)}_{\nu_t}\Big( kdW^n_n(t) - 2\sum_{i=1}^{n-1}dW^i_i(t)\Big) \nonumber \\
       &   & \quad +\ \sum_{i,j}S^{(k)}_{\nu_t}(\alpha_{\vec{l}}(t),\alpha_{\vec{m}}(t))
\frac{\langle\tau^j_i\alpha^{\vec{l}}(t), \alpha^{\vec{m}}(t)\rangle +\langle\alpha^{\vec{l}}(t),\tau^j_i\alpha^{\vec{m}}(t)\rangle}{\|\alpha(t)\|^{2}}
dW^i_j(t)  \nonumber \\
       &   & \quad +\ \frac{(n+1)k(n-k)\mu_2}{2n(n+2)}\Tr S^{(k)}_{\nu_t}dt,
\end{eqnarray}
where the Einstein summation convention is carried over the
indices $\vec{l}$ and $\vec{m}$.

We now have everything we need to present the main result of this
paper.

\begin{theorem}
\label{theorem.main}
Retain the assumptions and notation of Lemma \ref{theorem.ito.norm}.
Let ${\cal L}_k$ be the Lipschitz-Killing curvatures defined by \eqref{equation.LK}.
Then the It\^{o} derivatives of the Lipschitz-Killing curvatures for
$1\le k\le n-1$ are given by

\begin{eqnarray}
\label{equation.main}
d{\cal L}_{n-k-1}(M_t)& = &
\Big[K_{n,k}\int_{M}\int_{S(\mathbb{R})}\Big(
\sum_{i,\;p=1}^k S^{(k-1)}_{\nu_t}(\alpha_{\vec{l_p}}(t),\alpha_{\vec{m_i}}(t))\langle dB(u_{l_p}(t),u_{m_i}(t)),\nu_t\rangle \nonumber\\
                      &   & \mbox{} \times \langle\alpha^{\vec{l}}(t),\alpha^{\vec{m}}(t)\rangle \|\alpha(t)\|^{-1}\nonumber\\
                      &   & \mbox{} +\Tr S^{(k)}_{\nu_t}\|\alpha(t)\|\big( kdW^n_n -\sum_{i=1}^{n-1}dW^i_i\big)  \nonumber\\
              &   & \mbox{} +\sum_{i,\;j}S^{(k)}(\alpha_{\vec{l}},\alpha_{\vec{m}}) (\langle\tau^j_i\alpha^{\vec{l}}(t),\alpha^{\vec{m}}(t)\rangle +\langle\alpha^{\vec{l}}(t),\tau^j_i\alpha^{\vec{m}}(t)\rangle)dW^i_j \|\alpha(t)\|^{-1}\Big)\nonumber\\
              &   & \mbox{} \times 1_{N_{x_t}M_t}(-\nu_t){\cal H}_0(d\nu_t) {\cal H}_{n-1}(dx)\Big] \nonumber\\
              &   & \mbox{}+\frac{(n-k-1)(n+1)(k+1)\mu_2}{2n(n+2)}{\cal L}_{n-k-1}(M_t) dt,
\end{eqnarray}
where
${\cal H}_k$ is
 $k$-dimensional Hausdorff measure and
$N_{x_t}M_t$ is the normal cone to $M_t$ at $x_t\in M_t$.
\end{theorem}

Before giving the proof, we note that \eqref{equation.main} simplifies
considerably when $k=0$,  in which case, as we have already noted,
 $\lips_{n-1}(M_t)\equiv {\cal H}_{n-1}(M_t)$.
Then
\begin{eqnarray*}
{\cal L}_{n-1}(M_t) & = & K(n,0)\int_M\int_{S(\mathbb{R})}\|\alpha(t)\| 1_{N_{x_t}M_t}(-\nu_t){\cal H}_0(d\nu_t) {\cal H}_{n-1}(dx) \\
                    & = & \int_M \|\alpha(t)\| {\cal H}_{n-1}(dx).
\end{eqnarray*}
Hence, by Theorem~\ref{theorem.ito.norm},
\begin{eqnarray*}
d{\cal L}_{n-1}(M_t) & = & \int_M \|\alpha(t)\|\Big(\sum_{i=1}^{n-1}dW^i_i(t)
\Big) {\cal H}_{n-1}(dx) + \frac{(n-1)(n+1)\mu_2}{2n(n+2)}{\cal L}_{n-1}(M_t)dt.
\end{eqnarray*}

{\it En passant} to the proof we remind the reader of the dependence
of the various integrands appearing in \eqref{equation.main} on the
space parameter $x\in M$ through the vector field $U$, its various
spatial derivatives $W$ and $B$, and the tangent vectors, at $x\in
M$.

\textbf{Proof of Theorem \ref{theorem.main}:} We start with
\begin{eqnarray*}
d(\Tr S^{(k)}_{\nu_t}\|\alpha(t)\|) & = & (d(\Tr S^{(k)}_{\nu_t})\|\alpha(t)\| +\Tr S^{(k)}_{\nu_t}d(\|\alpha(t)\|) + \langle d\Tr S^{(k)}_{\nu_t},d\|\alpha(t)\|\rangle\\
                                  &  \definedas& I + II + III.
\end{eqnarray*}
We shall obtain a closed form expression for each of these terms.
By \eqref{equation.ito.trace} we have

\begin{eqnarray}
\label{equation.I}
I & = & (d(\Tr S^{(k)}_{\nu_t})\|\alpha(t)\| \nonumber \\
  & = & \Big[ \sum_{i,p=1}^k S^{(k-1)}_{\nu_t}(\alpha_{\vec{l_p}}(t),\alpha_{\vec{m_i}}(t))\langle dB(u_{l_p}(t),u_{m_i}(t)),\nu_t\rangle\Big] \langle\alpha^{\vec{l}}(t),\alpha^{\vec{m}}(t)\rangle\|\alpha(t)\|^{-1} \nonumber\\
           &   & \mbox{} +\Tr S^{(k)}_{\nu_t}\Big[ kdW^n_n(t) - 2\sum_{i=1}^{n-1}dW^i_i(t)\Big]  \|\alpha(t)\| \nonumber \\
       &   & \mbox{} +\sum_{i,j}S^{(k)}_{\nu_t}(\alpha_{\vec{l}}(t),\alpha_{\vec{m}}(t))(\langle\tau^j_i\alpha^{\vec{l}}(t), \alpha^{\vec{m}}(t)\rangle +\langle\alpha^{\vec{l}}(t),\tau^j_i\alpha^{\vec{m}}(t)\rangle)dW^i_j(t) \|\alpha(t)\|^{-1} \nonumber \\
       &   & \mbox{} +\frac{(n+1)k(n-k)\mu_2}{2n(n+2)}\Tr S^{(k)}_{\nu_t}\|\alpha(t)\|dt.
\end{eqnarray}

Using Theorem~\ref{theorem.ito.norm} we find
\begin{eqnarray}
\label{equation.II}
II & = & \Tr S^{(k)}_{\nu_t}d(\|\alpha(t)\|) \nonumber \\
   & = & \Tr S^{(k)}_{\nu_t}\|\alpha(t)\|\Big(\sum_{i=1}^{n-1}dW^i_i(t) +\frac{(n-1)(n+1)\mu_2}{2n(n+2)}dt\Big).
\end{eqnarray}

Finally using \eqref{equation.quad.cov}, \eqref{equation.ito.trace} and Theorem~\ref{theorem.ito.norm} we have
\begin{eqnarray}
\label{equation.III}
&&III \\
&&= \left\langle d\Tr S^{(k)}_{\nu_t},d\|\alpha(t)\|\right\rangle \nonumber\\
    &&= \Tr S^{(k)}_{\nu_t}\|\alpha(t)\|\Big\langle\Big( kdW^n_n -2\sum_{i-1}^{n-1}dW^i_i\Big) ,\, \sum_{i=1}^{n-1}dW^i_i\Big\rangle_t \nonumber\\
    &&  +\sum_{i,j}S^{(k)}_{\nu_t}(\alpha_{\vec{l}}(t),\alpha_{\vec{m}}(t))\big(\langle\tau^j_i\alpha^{\vec{l}}(t), \alpha^{\vec{m}}(t)\rangle +\langle\alpha^{\vec{l}}(t),\tau^j_i\alpha^{\vec{m}}(t)\rangle\big)
\Big\langle dW^i_j, \sum_{i=1}^{n-1}dW^i_i\Big\rangle\|\alpha(t)\|^{-1} \nonumber \\
    &&= \mbox{} -\frac{(n-1)(k+2)\mu_2}{n(n+2)}\Tr S^{(k)}_{\nu_t}\|\alpha(t)\|dt +\frac{2(n-k-1)\mu_2}{n(n+2)}\Tr S^{(k)}_{\nu_t}\|\alpha(t)\|dt.
\end{eqnarray}

Summing \eqref{equation.I}, \eqref{equation.II} and
\eqref{equation.III} we have
\begin{eqnarray*}
&&d(\Tr S^{(k)}_{\nu_t}\|\alpha(t)\|) \\
&&\qquad =\  \Big[ \sum_{i,p} S^{(k-1)}_{\nu_t}(\alpha_{\vec{l_p}}(t),\alpha_{\vec{m_i}}(t))\langle dB(u_{l_p}(t),u_{m_i}(t)),\nu_t\rangle\Big] \langle\alpha^{\vec{l}}(t),\alpha^{\vec{m}}(t)\rangle\|\alpha(t)\|^{-1} \\
&&\qquad\quad  + \ \Tr S^{(k)}_{\nu_t}\|\alpha(t)\|\Big[ kdW^n_n(t) - \sum_{i=1}^{n-1}dW^i_i(t)\Big]  \\
 &&\qquad\quad  + \ \sum_{i,j}S^{(k)}_{\nu_t}(\alpha_{\vec{l}}(t),\alpha_{\vec{m}}(t))\big(\langle\tau^j_i\alpha^{\vec{l}}(t), \alpha^{\vec{m}}(t)\rangle +\langle\alpha^{\vec{l}}(t),\tau^j_i\alpha^{\vec{m}}(t)\rangle\big)dW^i_j(t) \|\alpha(t)\|^{-1} \\
 &&\qquad\quad  +\ \frac{(n-k-1)(n+1)(k+1)\mu_2}{2n(n+2)} \Tr S^{(k)}_{\nu_t}\|\alpha(t)\|dt.
\end{eqnarray*}

Substituting the above  in \eqref{equation.LK} gives
\eqref{equation.main}
and so the theorem.  \qed \\

We can now easily deduce Theorem~\ref{theorem.mean.lk} of the Introduction.

\textbf{Proof of Theorem~\ref{theorem.mean.lk}:} In
\eqref{equation.main}, we note that with the single exception of the
last term, all terms are zero mean martingales due to the presence
of the martingale integrators $dW_.^.(t)$ or $dB^._{..}(t)$.
Therefore, taking expectations in
 \eqref{equation.main}, after taking the integral over time $t$,
immediatley yields
\beqq
\E \left\{ {\cal L}_{n-k-1}(M_t)
\right\}\ =\  \frac{(n-k-1)(n+1)(k+1)\mu_2}{2n(n+2)}
\int_0^t \E \left\{{\cal L}_{n-k-1}(M_s)\right\}\,ds.
\eeqq
Solving this linear differential equation gives
\eqref{equation.meanLK}, and we are done. \qed

We close with one further result, and two  open problems.

As for the first open problem, we believe that, in the setting of
Theorem~\ref{theorem.mean.lk}, \beqq \lim_{t\to\infty}
\frac{1}{t}\log \left(\frac{{\cal L}_{n-k-1}(M_t)}{{\cal
L}_{n-k-1}(M)}\right) \ = \ \frac{(n-k-1)(n+1)(k+1)\mu_2}{2n(n+2)},
\eeqq where the limit here is in $L_1$. At this point we do not have
an air tight proof of this. We would also like to add that the
almost sure growth rates may well be different from the ones
conjectured above.

As for the additional result, recall that throughout the paper, we
have assumed that $M$ was a codimension one manifold in $\real^n$.
From a technical point of view, this has a substantial simplifying
effect on the definition \eqref{LK:equation} of \LKC s. If
$\dim(M)=m<(n-1)$, then \eqref{LK:equation} changes in that the
normal cones are now of dimension $(n-m)$, $S(\real)$ is replaced by
$S(\real^{n-m})$ and so is of dimension $(n-m-1)$, and ${\cal H}_0$
and ${\cal H}_{n-1}$ are replaced by ${\cal H}_{n-m-1}$ and ${\cal
H}_{m}$, respectively. (The constant also changes, but this is less
important. See \cite{GRF} for details.) All told, we have \beq
\label{mk:equation} {\cal L}_{m-k}(M) \ =\
 K_{m,k}'\int_{M}\int_{S(\mathbb{R}^{m-m})}\Tr S^{(k)}_{\nu}
 1_{N_xM}(-\nu)\, {\cal H}_{n-m-1}(d\nu){\cal H}_{m}(dx),
\eeq
for some constants $ K_{m,k}'$. For all $k\neq 0$ here, we have found that
the complications introduced by increasing the dimension of the normal
space are such that computations analogous to those we have carried out
are forbiddingly complex.

For $k=0$ the  trace term in \eqref{mk:equation} disappears, and so it is
not hard to show that
\beq
d{\cal L}_{m}(M_t) & = & \int_M \|\alpha(t)\|\Big(\sum_{i=1}^{m}dW^i_i(t)\Big)
 {\cal H}_{m}(dx) + \frac{m(n-m)(n+1)\mu_2}{2n(n+2)}{\cal L}_{m}(M_t)dt.
\eeq
Since ${\cal L}_{m}(M_t)$ is the $m$-dimensional content of $M_t$, this is
a far from uninteresting result.

Of course, as before, this implies that
\beq
\E\left\{{\cal L}_{m}(M_t)\right\} \ =\
 {\cal L}_{m}(M)\exp\left(\frac{m(n-m)(n+1)\mu_2t}{2n(n+2)}\right).
\eeq
We have not, however, been able to find corresponding results for the
more general case.

{\smc Sreekar Vadlamani:} Faculty of Industrial Engineering
and Management, Technion - Israel Institute of Technology, Haifa, Israel
\newline sreekar@ieadler.technion.ac.il, http://tx.technion.ac.il/$\sim$ sreekar/

{\smc Robert J.\ Adler:} Faculty of Industrial Engineering
and Management, Technion - Israel Institute of Technology, Haifa, Israel
\newline robert@ieadler.technion.ac.il,
http://ie.technion.ac.il/Adler.phtml

\end{document}